\input amstex.tex
\magnification=1200
\documentstyle{amsppt}
\vcorrection{-1.0cm}
\topmatter
\title
     Analysis on the minimal representation of $O(p,q)$
\\
     --  I. Realization via conformal geometry
\endtitle
\author
  Toshiyuki KOBAYASHI
  and
  Bent \O RSTED
\endauthor
\affil
RIMS Kyoto and SDU-Odense University
\endaffil 
\address
 RIMS, Kyoto University,
   Sakyo-ku, Kyoto, 606-8502, Japan
\newline\indent
  Department of Mathematics and Computer Science,
  SDU - Odense University, Campusvej 55, DK-5230, Odense M, Denmark
\endaddress
\email{toshi\@kurims.kyoto-u.ac.jp;
       orsted\@imada.sdu.dk}\endemail
\abstract
This is the first in a series of papers devoted to
an analogue of the metaplectic representation,
namely the minimal unitary representation
of an indefinite orthogonal group;
this representation corresponds to the minimal
nilpotent coadjoint orbit
 in the philosophy of Kirillov-Kostant.
We begin by applying methods from conformal geometry
of pseudo-Riemannian manifolds to  a general construction
of an infinite-dimensional representation of the conformal
group on the solution space of the Yamabe equation.
By functoriality of the constructions, we obtain different
models of the unitary representation, as well as
giving new proofs of unitarity and irreducibility. 
The results in this paper play a basic role in the
subsequent papers, where we give explicit branching
formulae, and prove unitarization in the various 
models.  
\endabstract
\endtopmatter
\NoRunningHeads
\overfullrule=0pt

\def \Conf{\operatorname{Conf}}
\def \Isom{\operatorname{Isom}}

\def \Hsum#1{{{\underset{#1}\to{{\sum}^{\oplus}}}}}
\def \Ker{\operatorname{Ker}}

\def \trans{{}^t \!}
\def \tilLap#1{{\widetilde {\Delta}}_{#1}}
\def\fLap{\square_{\Bbb R^{p-1,q-1}}}

\def\Ad{\operatorname{Ad}}
\def\diag{\operatorname{diag}}
\define \F#1#2{F({#1},{#2})}               %
\define \sgn{\operatorname{sgn}}

\def \Vpq{V^{p,q}}

\def\pip#1#2#3{\pi^{{#1},{#2}}_{+, {#3}}}
\def\pim#1#2#3{\pi^{{#1},{#2}}_{-, {#3}}}

\define \set#1#2{\{{#1}:{#2}\}}
\def \sec#1{{\vskip 0.5pc\noindent$\underline{{\bold{\ch.\sc.}} \; \text{{#1}}}$\enspace}}
\def \num{{\ch.\sc}} %
\def\rarrowsim{\smash{\mathop{\,\rightarrow\,}\limits
  ^{\lower1.5pt\hbox{$\scriptstyle\sim$}}}}
\def\Pmax{P^{\roman{max}}}
\def\Mmax{M^{\roman{max}}}
\def\Amax{A^{\roman{max}}}
\def\Nmax{N^{\roman{max}}}

\def\prince#1#2#3{{#1}\text{-}\operatorname{Ind}_{\Pmax}^G({#3} \otimes \Bbb C_{#2})}
\def\princeK#1#2{{\operatorname{Ind}_{\Pmax}^G({#2} \otimes \Bbb C_{#1})}}

\def\spr#1#2{{\Cal H^{#1}(\Bbb R^{#2})}}
\def \Dom{\operatorname{Dom}}

\define \zdfc_#1^#2#3{{\Cal R}_{\frak#1}^{#2} ({\Bbb C}_{#3})}
\define \zdf_#1^#2{{\Cal R}_{\frak {#1}}^{#2}}   %

\define \rootsys#1#2{\triangle(#1,#2)}

\define \Ato#1#2{{\Cal A}\left({#1}\triangleright{#2}\right)}

\define \possys#1#2{\triangle^+(#1,#2)}
\def \Adisc{{A'}}

\def \gk{$(\frak g, K)$}
\def \xbec{\cite{1}}
\def \xbz{\cite{2}}
\def \xbrkosI{\cite{3}}
\def \xbrkosII{\cite{4}}
\def \xbrkosIII{\cite{5}}
\def \xerdHigI{\cite{6}}
\def \xerdIntII{\cite{7}}
\def \xgs{\cite{8}}
\def \xhela{\cite{9}} %
\def \xhowe{\cite{10}}
\def \xhowetan{\cite{11}}
\def \xhuzhu{\cite{12}}
\def \xkobast{\cite{13}}
\def \xkupq{\cite{14}}
\def \xkdecomp{\cite{15}} %
\def \xkdecoalg{\cite{16}} %
\def \xkdecoass{\cite{17}} %
\def \xkmfjp{\cite{18}}
\def \xkmf{\cite{19}}
\def \xkdecoaspm{\cite{20}}
\def \xkhcrrest{\cite{21}}
\def \xkoII{\cite{23}}
\def \xkoIII{\cite{24}}
\def \xkohcrcras{\cite{22}}
\def \xkos{\cite{25}}
\def \xlee{\cite{26}}
\def \xsab{\cite{27}}
\def \xschlap{\cite{28}}
\def \xschm{\cite{29}}%
\def \xschmid{\cite{}}%
\def \xtoramin{\cite{30}}
\def        \xvg{\cite{31}}%
\def \xvu{\cite{32}}
\def        \xvr{\cite{33}}%
\def \xvsing{\cite{34}}
\def \xvi{\cite{35}}
\def      \xwong{\cite{36}}%
\def \xorslmp{\cite{37}}
\def \xorsjfa{\cite{38}}

\document

\centerline{{\bf {Contents}}}
\smallskip
\item{ \S 1.}\enspace
 Introduction
\item{ \S 2.}\enspace
 Conformal geometry
\item{ \S 3.}\enspace
 Minimal unipotent representations of $O(p,q)$

\noindent
\def \ch{1}
\head
 \S \ch. Introduction
\endhead

\def \sc{0}
\sec{}
This is the first in a series of papers devoted to a study
of the so-called minimal representation of the semisimple Lie
group $G = O(p,q)$. We have taken the point of view that a rather
complete treatment of this representation and its various
realizations can be done in a self-contained way; also, such
a study involves many different tools from other parts of mathematics, 
such as differential geometry (conformal geometry and
pseudo-Riemannian geometry), analysis of solution spaces of
ultrahyperbolic differential equations, Sobolev spaces, special functions
such as hypergeometric functions of two variables,  
Bessel functions, 
analysis on semisimple symmetric spaces, 
and Dolbeault cohomology groups. Furthermore, the representation theory
yields new results back to these areas, so we feel it is worthwhile 
to illustrate such an interaction in as elementary a way as possible. 
The sequel (Part II) to the present paper contains sections 4-9, and
we shall also refer to these here. 
Part III is of more independent nature.

Working on a single unitary representation we essentially want to
analyze it by understanding its restrictions to natural subgroups, 
and to calculate intertwining operators between the various models -
all done very explicitly.
 We are in a sense studying
the symmetries of the representation space by breaking the large
symmetry present originally with the group $G$ by passing to a
subgroup. Geometrically the restriction is from the conformal
group $G$ to the subgroup of isometries $H$, where different
geometries (all locally conformally equivalent) correspond to
different choices of $H$. Changing $H$ will give rise to
radically different models of the representation, and at the
same time allow calculating the spectrum of $H$.      
       
Thus the overall aim is to elucidate as many aspects as possible of
a distinguished unitary irreducible representation of $O(p,q)$,
including its explicit branching laws to
natural subgroups and its explicit inner product on each geometric model.
 Our approach is also useful in understanding the
relation between the representation and a certain coadjoint
orbit, namely the minimal one, in the dual of the Lie algebra.
In order to give a good view of the perspective in our papers,
we are giving below a rather careful introduction to all these
aspects.

For a semisimple Lie group $G$ a particularly interesting unitary irreducible
representation, sometimes
called the minimal representation,
is the one corresponding via \lq\lq{geometric quantization}\rq\rq\ 
to the minimal nilpotent coadjoint orbit. It is still
a little mysterious in the present status of
the classification problem
 of the unitary dual of semisimple Lie groups. In recent years
several authors have considered the minimal
representation, and provided many new results,
in particular, Kostant, Torasso, Brylinski,
Li, Binegar, Zierau, and Sahi, 
 mostly by algebraic methods. 
For the double cover of the
symplectic group, this is the metaplectic representation, introduced
many years ago by Segal, Shale, and Weil. 
The explicit treatment
 of the metaplectic representation
 requires various methods from analysis and geometry, 
in addition to the
algebraic methods; and it is our aim in a series of papers to present for the
case of $G = O(p,q)$ the aspects pertaining to branching laws.  
{}From an algebraic view point of representation theory,
 our representations $\varpi^{p,q}$ are:
\newline\indent
i) minimal representations if $p + q \ge 8$ (i.e. the annihilator is the Joseph ideal).
\newline\indent 
ii) {\bf not} spherical if $p\neq q$ (i.e. no non-zero $K$-fixed vector).
\newline\indent
iii) {\bf not} highest weight modules of $SO_0(p,q)$ if $p, q \ge 3$.

Apparently our case
provides examples of new phenomena in representation theory, and we think 
that several aspects of our study can be applied to other cases as well.     
The metaplectic representation has had many applications
in representation theory and in number theory. A particularly
useful concept has been Howe's idea of dual pairs, where one
considers a mutually centralizing pair of subgroups in the
metaplectic group and the corresponding restriction of the metaplectic
representation. In Part II of our papers,
 we shall initiate a similar study
of explicit branching laws for other groups and representations
analogous to the classical case of Howe. 
Several such new examples of dual pairs have been
studied in recent years, mainly by algebraic techniques. Our case of 
the real orthogonal group presents a combination of abstract
representation theory and concrete analysis using methods from
conformal differential geometry. Thus we can relate the branching
law to a study of the Yamabe operator and its spectrum in
locally conformally equivalent manifolds; furthermore, we can
prove the existence of and construct explicitly  
an infinite discrete spectrum in the case where both factors in the dual pair
are non-compact.

The methods we use are further motivated by the theory of spherical harmonics,
extending analysis on the sphere to analysis on hyperboloids, and at the
same time using elliptic methods in the sense of analysis on
complex quadrics and the
theory of Zuckerman-Vogan's derived functor modules
 and their Dolbeault cohomological
realizations. 
Also important are general results on discrete
decomposability of representations and explicit knowledge of
branching laws.

It is noteworthy, that as we have indicated, this representation and its theory
of generalized Howe correspondence, illustrates several interesting
aspects of modern representation theory. Thus we have tried to be
rather complete in our treatment of the various models of the
representations occurring in the branching law. See for example
Fact 5.4, where we give three realizations: derived functor modules or
Dolbeault cohomology groups, eigenspaces on semisimple symmetric spaces,
and quotients of generalized principal series,
of the representations attached to minimal elliptic orbits. 

Most of the results of Part I and Part II were announced
 in \xkohcrcras, and
the branching law in the discretely decomposable case
(Theorem 7.1) was obtained in 1991, from which our study grew out. We have
here given the proofs of the branching laws for the minimal unipotent
representation and postpone the detailed treatment of the corresponding
classical orbit picture as announced in \xkohcrcras\ 
 to another paper. Also, the branching laws for
the representations associated to minimal elliptic orbits will appear
in another paper by one of the authors.

It is possible that a part of our results could be obtained
by using sophisticated results from the theory of dual pairs
in the symplectic group, for example the see-saw
rule (for which one may let our representation
correspond to the trivial representation of
one $SL(2, \Bbb R)$ member of the dual pair (\xhowetan, \xhuzhu).
We emphasize however, that our approach is quite explicit and has
the following advantages:

 (a) It is not only an abstract representation theory
 but also attempts new interaction
 of the minimal representation with analysis
 on manifolds.  
For example,
in Part II we use in an elementary way
conformal differential geometry and the functorial properties
of the Yamabe operator to construct the minimal representation
and the branching law in a way which seems promising for other
cases as well;
each irreducible constituent is explicitly constructed by using
explicit intertwining operators via local conformal diffeomorphisms
between spheres and hyperboloids.

(b) For the explicit intertwining operators we obtain
Parseval-Plancherel type theorems, i.e. explicit $L^2$ versions
of the branching law and the generalized Howe correspondence.  
This also  gives a good
perspective on the continuous spectrum, in particular
yielding a natural conjecture for the complete Plancherel
formula.

A special case of our branching law illustrates the physical situation
of the conformal group of space-time $O(2,q)$; here the minimal
representation may be interpreted either as the mass-zero 
spin-zero wave equation, or as the bound states of the Hydrogen atom
(in $q-1$ space dimensions). Studying the branching law means
breaking the symmetry by for example restricting to the isometry
group of De Sitter space $O(2,q-1)$ or anti De Sitter space $O(1,q)$.
In this way the original system (particle) breaks up into
constituents with less symmetry.

In Part III,
 we shall realize the same representation on a
 space of solutions of the ultrahyperbolic equation
 $\fLap f =0$ on $\Bbb R^{p-1,q-1}$,
 and give an intrinsic inner product as an
 integration over a non-characteristic hypersurface.

Completing our discussion of different models of the
minimal representation,
 we find yet another explicit intertwining operator,
this time to an $L^2$ - space of functions on a hypersurface (a cone)
in the nilradical of a maximal parabolic $P$ in $G$. We find the
$K$-finite functions in the case of $p+q$ even in terms of
modified Bessel functions.
We remark that Vogan pointed out
 a long time ago
 that there is no minimal representation
 of $O(p,q)$ if $p+q >8$ is odd \xvsing. 
 On the other hand,
  we have found a new interesting phenomenon that
 in the case $p+q$ is odd 
 there still exists a geometric model of 
 a \lq\lq minimal representation\rq\rq\
 of $\frak {o}(p,q)$ with a natural inner product (see Part III).
 Of course, such a representation 
 does not have non-zero $K$-finite vectors  for $p+q$ odd,
 but have $K'$-finite vectors for smaller $K'$.  
 What we construct in this case is an element of the category of 
 $({\frak g},P)$ modules
 in the sense that it globalizes to $P$ (but not $K$); 
 we feel this concept perhaps plays
a role for other cases of the orbit method as well. 
\par
In summary, we give a geometric and intrinsic model of the
minimal representation 
$\varpi^{p,q}$
(not coming from the construction of $\varpi^{p,q}$ by the
$\theta$-correspondence)
on $S^{p-1} \times S^{q-1}$ and on various 
pseudo-Riemannian manifolds which are conformally equivalent, 
using the functorial properties of the Yamabe operator,
a key element in conformal differential geometry. The
branching law for 
$\varpi^{p,q}$ gives at the same time new perspectives 
in conformal geometry, and we give an
independent proof of the unitarity  
of a subspace of
the kernel of the Yamabe operator.
The main interest in this special
case of a small unitary representation is not only 
to obtain the formulae, but also to investigate the
geometric
and analytic methods, which provide new ideas in representation
theory. 
  
Leaving the general remarks, let us now for the rest of this
introduction be a little more specific about the
contents of the present paper.

\def \sc{1}
\sec{}
Let $G$ be a reductive Lie group,
 and $G'$ a reductive subgroup of $G$.
We denote by $\widehat{G}$  the unitary dual of $G$,
 the equivalence classes of irreducible unitary representations of $G$.
Likewise $\widehat{G'}$ for $G'$.
If $\pi \in \widehat{G}$,
 then the restriction $\pi|_{G'}$ is not necessarily irreducible.
By a branching law,
 we mean an explicit irreducible decomposition formula:
$$
   \pi|_{G'} \simeq \int_{\widehat{G'}}^\oplus m_\pi(\tau) \tau d \mu(\tau)
   \qquad
   \text{(direct integral)},
\tag \num.1
$$
 where $m_\pi(\tau) \in \Bbb N \cup \{\infty\}$ and $d \mu$ is a
 Borel measure on $\widehat{G'}$.
\def \sc{2}
\sec{}
We denote by $\frak g_0$ the Lie algebra of $G$.
The {\bf orbit method} due to Kirillov-Kostant
 in the unitary representation theory of Lie groups
 indicates that  
 the coadjoint representation
 $\Ad^* \:G  \to GL(\frak g_0^*)$ 
 often has a surprising intimate relation
 with the unitary dual $\widehat G$.
It works perfectly for simply connected nilpotent Lie groups.
For real reductive Lie groups $G$,
 known examples suggest that
  the set of coadjoint orbits $\sqrt{-1}\frak g_0^*/G$
  (with certain integral conditions)
 still gives a fairly good approximation 
 of the unitary dual $\widehat G$.

\def \sc{3}
\sec{}
Here is a rough sketch of a unitary representation $\pi_\lambda$ of $G$,
 attached to an elliptic element $\lambda \in \sqrt{-1} \frak g_0^*$:
The elliptic coadjoint orbit $\Cal O_\lambda = \Ad^*(G) \lambda$ 
 carries a $G$-invariant complex structure,
 and one can define a $G$-equivariant holomorphic line bundle
 $\widetilde{\Cal L_\lambda} :=
 \Cal L_\lambda \otimes (\wedge^{top} T^* \Cal O_\lambda)^\frac12$ over $\Cal O_\lambda$,
 if $\lambda$ satisfies some integral condition.
Then,
 we have a Fr\'echet representation of $G$
 on the Dolbeault cohomology group
 $H_{\bar\partial}^S(\Cal O_\lambda, \widetilde{\Cal L_\lambda})$,
 where $S := \dim_\Bbb C \Ad^*(K) \lambda$ (see \xwong\ for details),
 and of which a unique dense subspace we can define
 a unitary representation $\pi_\lambda$ of $G$ (\xvu)
 if $\lambda$ satisfies certain positivity.
The unitary representation $\pi_\lambda$ is irreducible and non-zero
 if $\lambda$ is sufficiently regular.
The underlying $(\frak g, K)$-module is so called
 \lq\lq $A_{\frak q}(\lambda)$\rq\rq\ in the sense of Zuckerman-Vogan
 after certain $\rho$-shift.

In general,
 the decomposition (\ch.1.1) contains both discrete and continuous spectrum.
The condition for the discrete decomposition (without continuous spectrum)
 has been studied in  \xkdecomp, \xkdecoalg, \xkdecoass, and \xkdecoaspm,
 especially for $\pi_\lambda$ attached to elliptic orbits $\Cal O_\lambda$.
It is likely that if $\pi \in \widehat{G}$ is \lq\lq attached to\rq\rq\
 a nilpotent orbit,
 which is contained in the limit set of $\Cal O_\lambda$,
 then the discrete decomposability of $\pi|_{G'}$
 should be inherited from that of the elliptic case $\pi_\lambda|_{G'}$.
We shall see in Theorem~4.2 
 that this is the case in our situation.
\def \sc{4}
\sec{}
There have been a number of attempts to construct
 representations attached to nilpotent orbits.
Among all,
 the Segal-Shale-Weil representation (or the oscillator representation)
 of $\widetilde{Sp}(n,\Bbb R)$,
 for which we write $\tilde\varpi$,
 has been best studied,
 which is supposed to be attached to the minimal nilpotent orbit
 of $\frak {sp}(n,\Bbb R)$.
The restriction of $\tilde\varpi$ to a reductive dual pair $G' = G_1' G_2'$
 gives Howe's correspondence (\xhowe).

The group $\widetilde{Sp}(n,\Bbb R)$ is a split group of type $C_n$,
 and analogous to $\tilde\varpi$,
 Kostant constructed a minimal representation of $SO(n,n)$,
 a split group of type $D_n$.
Then Binegar-Zierau generalized it for $SO(p,q)$ with $p+q \in 2 \Bbb N$.
This representation (precisely, of $O(p,q)$, see Section 3) will be denoted by
 $\varpi^{p,q}$.

\def \sc{5}
\sec{}
Let $G' := G_1' G_2' = O(p',q') \times O(p'', q'')$,
 $(p' + p''=p, q' + q'' = q)$,
 be a subgroup of $G = O(p,q)$.
Our object of study in Part II will be the branching law ${\varpi^{p,q}}_{|{G'}}$.
We note that $G_1'$ and $G_2'$ form a  mutually centralizing pair of
 subgroups in $G$.

It is interesting to compare the feature of the following two cases:
\newline
(i)\enspace\
 the restriction $\tilde\varpi|_{G_1' G_2'}$
 (the Segal-Shale-Weil representation for type $C_n$),
\newline
(ii)\enspace
 the restriction $\varpi^{p,q}|_{G_1' G_2'}$
 (the Kostant-Binegar-Zierau representation for type $D_n$).

The reductive dual pair
$
   (G, G') = (G, G_1' G_2')
$
 is of the $\otimes$-type in (i),
 that is, induced from $GL(V) \times GL(W) \to GL(V \otimes W)$;
 is of the $\oplus$-type in (ii),
 that is, induced from $GL(V) \times GL(W) \to GL(V \oplus W)$.
On the other hand,
 both of the restrictions in (i) and (ii)
 are discretely decomposable if one factor $G_2'$ is compact.
On the other hand,
 $\tilde\varpi$ is (essentially) a highest weight module in (i),
 while $\varpi^{p,q}$ is not if $p, q > 2$ in (ii).

\def \sc{6}
\sec{}
Let  $p+q \in 2 \Bbb N$, $p, q \ge 2$, and $(p,q) \neq (2,2)$.
In this section we state the main results of the present paper and 
 the sequels (mainly Part II;
 an introduction of Part III will be given separately in \xkoIII).
 The first Theorem A below (Theorem 2.5) says that there
is a general way of constructing representations of a conformal
group by twisted pull-backs (see Section 2 for notation). It is the
main tool to give different models of our representation.  

\proclaim{Theorem~A}
Suppose 
 that a group  $G$ acts conformally on 
  a  pseudo-Riemannian manifold $M$ of dimension $n$.  
  \newline{\rm 1)}\enspace
  Then, 
  the Yamabe operator (see (2.2.1) for the definition)
  $$
          \tilLap{M} \: C^{\infty} (M) \to C^{\infty} (M)
          $$
           is an intertwining operator
                from $\varpi_{\frac {n-2} 2}$ to $\varpi_{\frac {n+2}2}$
         (see (2.5.1) for the definition of $\varpi_{\lambda}$).  
\newline{\rm 2)}\enspace
                The kernel $\Ker \tilLap{M}$ is a subrepresentation of $G$
                 through $\varpi_{\frac {n-2}2}$.  
\endproclaim

\proclaim{Theorem~B}
{\rm 1)}\enspace
The minimal representation $\varpi^{p,q}$ of $O(p,q)$
 is realized as the kernel of the Yamabe operator
 on $S^{p-1} \times S^{q-1}$.
\newline
{\rm 2)}\enspace 
 $\varpi^{p,q}$ is also realized as
 a subspace (roughly, half) of the kernel of
 the Yamabe operator on the hyperboloid 
 $\set{(x,y) \in \Bbb R^{p,q}}{|x|^2 - |y|^2 =1}$.
\newline
{\rm 3)}\enspace
$\varpi^{p,q}$ is also realized in a space
of solutions to the Yamabe equation on $\Bbb R^{p-1,q-1}$
which is a standard ultrahyperbolic constant coefficient
differential equation. 
\newline
{\rm 4)}\enspace
 $\varpi^{p,q}$ is also realized as the unique non-trivial subspace of
 the Dolbeault cohomology group
 $H_{\bar \partial}^{p-2} (G/L, \Cal L_{\frac{p+q-4}{2}})$
\endproclaim

In Theorem B~(1) is contained in Part I, Theorem~3.6.1;
 (2) in Part II, Corollary~7.2.1;
 and (3) in Part III, Theorem~4.7.
In each of these models,
 an explicit model is given explicitly.
 In the models (2) and (3),
  the situation is subtle because the \lq\lq action\rq\rq\ of 
  $O(p,q)$ is no more smooth but only meromorphic.
  Then Theorem~A does not hold in its original form,
   and we need to carry out a careful analysis for it 
   (see Part II and Part III).
The proof of the statement (4)
  will appear in another paper. 
Here $G/L$ is an elliptic coadjoint orbit as in
\S1.3, and $L = SO(2) \times O(p-2,q)$.

The branching laws in Theorem~C and Theorem~D are
the main themes in Part II; for notation see section 7 and section 9.

\proclaim{Theorem~C}
If $q'' \ge 1$ and $q' + q'' = q$,
 then the twisted pull-back $\widetilde{\Phi_1^*}$
 of the local conformal map $\Phi_1$
 between spheres
 and hyperboloids
 gives an explicit irreducible decomposition
 of the unitary representation $\varpi^{p,q}$
 when restricted to $O(p,q') \times O(q'')$:
$$
 \varpi^{p,q}|_{O(p, q') \times O(q'')}
 \simeq
 \overset{\infty\hphantom{M}}\to{\Hsum{l=0\hphantom{M}}}
  \pip{p}{q'}{l + \frac{q''}2 -1} \boxtimes \spr{l}{q''}.
$$
\endproclaim
In addition, we give in \S 8, Theorem~8.5 the Parseval-Plancherel theorem
for the situation in Theorem C on the \lq\lq hyperbolic space model\rq\rq. 
This may be also regarded as the unitarization of the minimal representation $\varpi^{p,q}$.

 The twisted pull-back for a locally conformal diffeomorphism
  is defined for an arbitrary pseudo-Riemannian manifold 
  (see Definition~2.3).

\proclaim{Theorem~D}
The twisted pull-back of the locally conformal diffeomorphism
 also constructs
$$
\Hsum{\lambda \in \Adisc(p', q') \cap \Adisc(q'',p'')}
  \pip{p'}{q'}{\lambda} \boxtimes \pim{p''}{q''}{\lambda}
\oplus
\Hsum{\lambda \in \Adisc(q', p') \cap \Adisc(p'',q'')}
  \pim{p'}{q'}{\lambda} \boxtimes \pip{p''}{q''}{\lambda}
$$
as a discrete spectra in the branching law.
\endproclaim 
\def \sc{7}
\sec{}
The papers (Part I and Part II) are organized as follows:
Section 2 provides a conformal construction of a representation on the kernel
 of a shifted Laplace-Beltrami operator.
In section 3, we construct an irreducible unitary representations,
 $\varpi^{p,q}$  of $O(p,q)$ ($p + q \in 2 \Bbb N, p, q \ge 2$)
\lq\lq attached to\rq\rq\
 the minimal nilpotent orbit
 applying Theorem~2.5.
This representation coincides with the minimal representation
 studied by Kostant, Binegar-Zierau (\xbz, \xkos).
In section 3 we give a new intrinsic characterization of the
Hilbert space for the minimal representation in this
model, namely as a certain Sobolev space of solutions, 
see Theorem 3.9 and Lemma 3.10. 
Such Sobolev estimate will be used in the construction of discrete spectrum
 of the branching law in section 9.
Section 4 contains some general results on discrete decomposable 
 restrictions (\xkdecoalg, \xkdecoass), specialized in
detail to the present case. Theorem 4.2 characterizes 
 which dual pairs in our situation provide discrete decomposable
  branching laws of the restriction of the minimal representation
  $\varpi^{p,q}$.
In section 5, we introduce unitary representations,
  $\pi_{\pm,\lambda}^{p,q}$ of $O(p,q)$ \lq\lq attached to\rq\rq\
  minimal elliptic coadjoint orbits.
In sections 7 and  9,
 we give a discrete spectrum of the branching law $\varpi^{p,q}|_{G'}$
 in terms of $\pi_{\pm, \lambda}^{p',q'} \in \widehat{O(p',q')}$
 and $\pi_{\pm, \lambda}^{p'',q''} \in \widehat{O(p'',q'')}$.
In particular,
 if one factor $G_2'=O(p'',q'')$ is compact (i.e\. $p''=0$ or $q''=0$),
 the branching law is completely determined 
 together with a Parseval-Plancherel theorem in section 8. 

Notation:
 $\Bbb N = \{ 0, 1, 2, \dots \}$.

The first author expresses his sincere gratitude to SDU - Odense University for
 the warm hospitality. 
\def \ch{2}
\def \sc{1}
\head
 \S \ch. Conformal geometry
\endhead
\sec{}
The aim of this section
 is to associate a distinguished
 representation $\varpi_M$ of the conformal group
 $\Conf (M)$ to a general pseudo-Riemannian manifold $M$
 (see Theorem \ch.5).

\def \sc{2}
\sec{}
Let $M$ be an $n$ dimensional manifold
 with pseudo-Riemannian metric $g_M$
 $(n \ge 2)$.
Let $\nabla$ be the Levi-Civita connection
 for the pseudo-Riemannian metric $g_M$.
The curvature tensor field $R$ is defined by
$$
    R(X,Y)Z:= \nabla_X \nabla_Y Z -\nabla_Y \nabla_X Z -\nabla_{[X,Y]} Z,
    \quad
    X, Y, Z \in \frak X(M).
$$
We take an orthonormal basis
 $\{X_1, \cdots, X_n\}$ of $T_x M$ for a fixed $x \in M$.
Then the scalar curvature $K_M$ is defined by
$$
     K_M(x):= \sum_{i=1}^n \sum_{j=1}^n g_M \left(R(X_i,X_j)X_i,X_j\right).
$$
The right side is independent of the choice of the basis $\{X_i\}$
 and so $K_M$ is a well-defined function on $M$.
We denote by $\Delta_M$
 the Laplace-Beltrami operator on $M$.
The {\bf Yamabe operator} is defined to be
$$
    \tilLap{M} := \Delta_M - \frac {n-2}{4(n-1)} K_M.
\tag \num.1
$$
See for example \xlee\ for a good discussion of the geometric meaning
and applications of this operator.
Our choice of the signature of $K_M$ and $\Delta_M$ is illustrated
 as follows:
\example{Example \num}
We equip $\Bbb R^n$ and $S^n$ with standard Riemannian metric.
Then
$$
\alignat3
     &\text{For } \Bbb R^n;  \qquad\qquad
     &K_{\Bbb R^n} &\equiv 0, \qquad
     &\tilLap{\Bbb R^n}   
       &=\Delta_{\Bbb R^n}
       = \sum_{i=1}^n \frac {\partial^2}{\partial x_i^2}.
\\
     &\text{For } S^n;     \qquad\qquad
     & K_{S^n} &\equiv (n-1)n, \qquad
     &\tilLap{S^n}&=\Delta_{S^n} - \frac 1 4 n(n-2).
\endalignat
$$
\endexample    
\def \sc{3}
\sec{}
Suppose $(M, g_M)$ and $(N, g_N)$
 are pseudo-Riemannian manifolds of dimension $n$.  
A local diffeomorphism $\Phi \: M \to N$ is called
 a {\it{conformal map}}
 if there exists a positive valued function $\Omega$ on $M$
 such that
$$
     \Phi ^* g_N = \Omega ^2 g_M.  
$$
$\Phi$ is isometry if and only if $\Omega \equiv 1$
 by definition.

We  denote the group of conformal transformation (respectively, isometry)
 of a pseudo-Riemannian manifold $(M, g_M)$ by
$$
\alignat2
  &\Conf(M)&&:=\set{\Phi \in \operatorname{Diffeo}(M)}
                {\Phi \: M \to M \text{ is conformal}},     
\\
  &\Isom(M)&&:=\set{\Phi \in \operatorname{Diffeo}(M)}
                {\Phi \: M \to M \text{ is isometry}}.     
\endalignat
$$
Clearly, $\Isom(M) \subset \Conf(M)$.

If $\Phi$ is conformal,
 then we have
 (e\.g\. \xorslmp; \xhela, Chapter II, Excer\. A\.5)
$$
  \Omega^{\frac {n+2} 2} (\Phi^* \tilLap{N} f)
  = \tilLap{M} (\Omega^{\frac {n-2}2} \Phi ^* f) 
  \tag {\num.1}
$$
 for any $f \in C^\infty(N)$.
We define a twisted pull-back
$$
     \Phi_{\lambda}^* \: C^{\infty} (N) \to C^{\infty}(M),\quad
                         f \mapsto \Omega^{\lambda}(\Phi^* f),
  \tag{\num.2}
$$
 for each fixed $\lambda \in \Bbb C$.  
Then the formula (\num.1) is rewritten as 
$$
     \Phi_{\frac {n+2}{2}}^* \tilLap{N} f
    =\tilLap{M} \Phi_{\frac {n-2}{2}}^* f.  
  \tag {\num.1$'$}
$$
The case when $\lambda = \frac{n-2}2$ is particularly important.
Thus, we write the twisted pull-back for $\lambda = \frac{n-2}2$ as follows:
\definition{Definition \num}
$
  \widetilde {\Phi^*} = \Phi_{\frac {n-2}{2}}^*
  \: C^{\infty}(N) \to C^{\infty}(M),\quad
     f \mapsto \Omega^{\frac {n-2} 2} (\Phi^* f).  
$
\enddefinition
Then the formula (\num.1) implies that
$$
     \tilLap{N} f =0
     \quad
     \text{ on }\Phi(M)
     \quad
     \text{if and only if}
     \quad
     \tilLap{M} (\widetilde {\Phi^*}f) =0
     \quad
     \text{ on }M
  \tag {\num.3}
$$
because $\Omega$ is nowhere vanishing.  
\par
If $n=2$,
 then
 $\tilLap{M} =\Delta_M$,
 $\tilLap{N}=\Delta_N$,
 and
 $\widetilde{\Phi^*}=\Phi^*$.  
Hence,
 (\num.3) implies a well-known fact
 in the two dimensional case
 that {\sl a conformal map $\Phi$ preserves harmonic functions},
 namely,
$$
     \text{$f$ is harmonic}
     \Leftrightarrow
     \text{$\Phi^* f$ is harmonic.  }
$$

\def \sc{4}
\sec{}
Let $G$ be a Lie group acting conformally
 on a pseudo-Riemannian manifold $(M, g_M)$.
We write the action of $h \in G$ on $M$
 as $L_h \: M \to M, x \mapsto L_h x$.  
By the definition of conformal transformations,
 there exists a positive valued function $\Omega(h,x)$
 $(h \in G, x \in M)$ such that
$$
     L_h^* g_M = \Omega(h,\cdot)^2 g_M
     \quad (h \in G).  
$$
Then we have
\proclaim{Lemma \num}
For $h_1, h_2 \in G$
 and $x \in M$, 
 we have
$$
     \Omega(h_1 h_2, x)= \Omega(h_1, L_{h_2} x) \ \Omega(h_2,x).  
$$
\endproclaim
\demo{Proof}
It follow from $L_{h_1 h_2} = L_{h_1} L_{h_2}$
 that
$$
     L_{h_1 h_2}^* g_M = L_{h_2}^* L_{h_1}^* g_M.  
$$
Therefore we have
$
     \Omega(h_1 h_2, \cdot)^2 \ g_M
    = L_{h_1 h_2}^* g_M 
= L_{h_2}^* \left(L_{h_1}^* g_M\right)
    =L_{h_2}^* \left(\Omega(h_1, \cdot)^2 \ g_M \right)
$
\linebreak  
$
    =\Omega (h_1, L_{h_2} \cdot)^2 \ \Omega (h_2, \cdot)^2 \ g_M.  
$
Since $\Omega$ is a positive valued function,
 we conclude that
$
     \Omega (h_1 h_2, x)= \Omega(h_1, L_{h_2} x) \ \Omega(h_2, x).  
$
\qed
\enddemo

\def \sc{5}
\sec{}
For each $\lambda \in \Bbb C$,
 we form a representation $\varpi_{\lambda} \equiv \varpi_{M, \lambda}$
 of the conformal group $G$ on $C^{\infty} (M)$ as follows:
$$
     \left(\varpi_{\lambda}(h^{-1}) f\right) (x)
      := \Omega(h,x)^{\lambda} f (L_h x),
     \quad
     (h \in G, f \in C^{\infty}(M), x \in M).  
\tag \num.1
$$
Then Lemma \ch.4 assures
 that 
$
     \varpi_{\lambda}(h_1) \ \varpi_{\lambda}(h_2)=\varpi_{\lambda}(h_1 h_2),
$
 namely,
 $\varpi_{\lambda}$ is a representation of $G$.

Denote by $d x$ the volume element on $M$ defined by
 the pseudo-Riemannian structure $g_M$.
Then we have
$$
    L_h^* (d x) = \Omega(h, x)^n d x
    \quad \text{ for }
    h \in G.
$$
Therefore,
 the map $f \mapsto f \, d x$ gives a $G$-intertwining operator
 from $(\varpi_n, C^\infty(M))$ into the space of distributions $\Cal D'(M)$
 on $M$.

Here is a construction of a representation of the group %
 of conformal diffeomorphisms of $M$.  
\proclaim{Theorem \num}
Suppose 
 that a group  $G$ acts conformally on 
 a  pseudo-Riemannian manifold $M$ of dimension $n$.  
Retain the notation before.  
\newline{\rm 1)}\enspace
Then, 
the Yamabe operator
$$
    \tilLap{M} \: C^{\infty} (M) \to C^{\infty} (M)
$$
 is an intertwining operator
 from $\varpi_{\frac {n-2} 2}$ to $\varpi_{\frac {n+2}2}$.  
\newline{\rm 2)}\enspace
The kernel $\Ker \tilLap{M}$ is a subrepresentation of $G$
 through $\varpi_{\frac {n-2}2}$.  
\endproclaim
\demo{Proof}
(1) is a restatement of the formula (\ch.3.1).
(2) follows immediately from (1).
\qed
\enddemo
The representation of $G$ on $\Ker \tilLap{M}$
 given in Theorem \num (2) will be denoted by $\varpi \equiv \varpi_M$.

\def \sc{6}
\sec{}
Here is a naturality of the
 representation of the conformal group $\Conf(M)$
 in Theorem \ch.5:
\proclaim{Proposition~\num}
Let $M$ and $N$ be pseudo-Riemannian manifolds of dimension $n$,
 and a local diffeomorphism $\Phi: M \to N$ be a conformal map.
Suppose that Lie groups $G'$ and $G$  act conformally on 
 $M$ and $N$, respectively.
The actions will be denoted by $L_M$ and $L_N$, respectively.
We assume that there is a homomorphism $i\: G' \to G$ such that
$$
     L_{N,i(h)} \circ \Phi = \Phi \circ L_{M,h}
     \quad
     (\text{ for any } h \in G').  
$$

We write conformal factors $\Omega_M$, $\Omega_N$ and $\Omega$ as follows:
$$
\alignat2
     &L_{M, h}^* g_M &&= \Omega_M(h, \cdot)^2 g_M \qquad (h \in G'),
\\
     &L_{N, h}^* g_N &&= \Omega_N(h, \cdot)^2 g_N \qquad (h \in G),
\\
     &\Phi^* g_N &&= \Omega^2 g_M.
\endalignat
$$
\item{\rm{1)}} For $x \in M$ and $h \in G'$,
 we have
$$
   \Omega(L_{M, h} x) \ \Omega_M(h, x)
   =
   \Omega(x) \   \Omega_N(i(h), \Phi(x)).
\tag \num.1
$$
\item{\rm {2)}}
 Let $\lambda \in \Bbb C$ and
 $\Phi_{\lambda}^* : C^\infty(N) \to C^\infty(M)$
 be the twisted pull-back defined in (\ch.3.2).  
Then $\Phi_{\lambda}^*$ respects 
 the $G$-representation 
 $(\varpi_{N, \lambda}, C^\infty(N))$
 and the $G'$-representation 
 $(\varpi_{M, \lambda}, C^\infty(M))$
 through $i\: G' \to G$.  
\item{\rm {3)}} $\widetilde {\Phi^*} = \Phi_{\frac {n-2}2}^* \: 
 C^{\infty}(N) \to C^{\infty} (M)$
 sends $\Ker \tilLap{N}$ into $\Ker \tilLap{M}$.  
In particular,
 we have a commutative diagram:
$$
  \CD
        \Ker \tilLap{N} @>{\widetilde{\Phi^*}}>> \Ker \tilLap{M}
\\
         @V{\varpi_N(i(h))}VV                             @VV{\varpi_M(h)}V
\\
        \Ker \tilLap{N} @>>{\widetilde{\Phi^*}}> \Ker \tilLap{M}
  \endCD
  \tag{\num.2}
$$
 for each $h \in G'$.   
\item{\rm {4)}} 
If $\Phi$ is a diffeomorphism onto $N$,
 then $(\Phi^{-1})_\lambda^*$ is the inverse of $\Phi_\lambda^*$
 for each $\lambda \in \Bbb C$.
In particular,
 $\widetilde{\Phi^*}$ is a bijection 
 between $\Ker \tilLap{N}$ and $\Ker \tilLap{M}$ 
 with inverse $\widetilde{(\Phi^{-1})^*}$.
\endproclaim
\demo{Proof}
1)\enspace
Because $L_{N, i(h)} \circ \Phi = \Phi \circ L_{M, h}$ for $h \in G'$,
 we have
$$
   (\Phi^* L_{N, i(h)}^* g_N) (x) =
   (L_{M, h}^* \Phi^* g_N) (x), 
  \quad
 \text{for }
   x \in M.
$$
Hence,
$$
 \Omega_N(i(h), \Phi(x))^2 \   \Omega(x)^2 g_M(x)
  =
   \Omega(L_{M, h} x)^2 \ \Omega_M(h, x)^2 g_M(x).
$$
Because all conformal factors are positive-valued functions,
 we have proved (\num.1).
\newline{2)}\enspace
We want to prove
$$
( \varpi_{M,\lambda}(h^{-1}) \Phi_{\lambda}^* f)(x)
 = (\Phi_{\lambda}^* \varpi_{N,\lambda}(i(h^{-1})) f)(x)
\tag \num.3
$$
 for any $x \in M$, $h \in G'$ and $\lambda \in \Bbb C$.  
In view of the definition, we have
$$
\align
   \text{the left side of (\num.3)}
 = &(\varpi_{M,\lambda}(h^{-1}) (\Omega^{\lambda} \ \Phi^* f))(x)
\\
 =& \Omega_M(h, x)^\lambda \ \Omega (L_{M,h} x)^{\lambda}
     (\Phi^* f)(L_{M,h} x)
\\
 =&\Omega(x)^{\lambda} \Omega_N(i(h),\Phi(x))^{\lambda} 
  f(\Phi \circ L_{M,h} x).  
\intertext{Here the last equality follows from (\num.1).  }
 \text{The right side of (\num.3)}
 =& (\Phi_{\lambda}^* \Omega_N(i(h), \cdot)^\lambda \ f(L_{N,h} \cdot))(x)
\\
 =& \Omega(x)^{\lambda}\ \Omega_N(i(h), \Phi(x))^\lambda 
    \ f(L_{N,i(h)} \circ \Phi(x)).  
\endalign
$$
Therefore,
 we have (\num.3),
 because 
$
 L_{N,i(h)} \circ \Phi = \Phi \circ L_{M,h}.  
$
\newline{3)}\enspace
If $f \in C^{\infty}(N)$ satisfies $\tilLap{N}f=0$,
 then $\tilLap{M}(\widetilde{\Phi^*}f)
=\Omega^{\frac {n+2}2}(\Phi^* \tilLap{N} f)=0$
 by (\ch.3.1).  
Hence $\widetilde {\Phi^*}(\Ker \tilLap{N}) \subset \Ker \tilLap{M}$.  
The commutativity of the diagram (\num.2) follows from (2) and Theorem~\ch.5~(2),
 if we put $\lambda = \frac {n-2}2$.  
\newline{4)}\enspace
Because $(\Phi^{-1})^* g_M = (\Omega \circ \Phi^{-1})^{-2} g_N$,
 the twisted pull-back $(\Phi^{-1})_\lambda^* F$ is given by 
 the following formula from definition (\ch.3.2):
$$
(\Phi^{-1})_\lambda^*
  \: 
    C^\infty(M) \to C^\infty(N), \
    F \mapsto
   (\Phi^{-1})_\lambda^* F 
 = (\Omega \circ \Phi^{-1})^{-\lambda} (F \circ \Phi^{-1}).
$$ 
Now the statement (4) follows immediately.
\qed
\enddemo

\def \ch{3}
\def \sc{1}
\head
 \S \ch. Minimal unipotent representations of $O(p,q)$
\endhead
\sec{}
In this section,
 we apply Theorem~2.5 to the specific setting where
 $M = S^{p-1} \times S^{q-1}$ is equipped with
 an indefinite Riemannian metric,
 and where
 the indefinite orthogonal group $G = O(p,q)$ acts conformally on $M$.
The resulting representation, 
 denoted by $\varpi^{p,q}$, 
 is non-zero, irreducible and unitary if $p+q \in 2 \Bbb N, p, q \ge 2$
 and if $(p,q) \neq (2,2)$.
This representation coincides with the one constructed
 by Kostant, Binegar-Zierau (\xbz, \xkos),
 which has the Gelfand-Kirillov dimension $p+q-3$ (see Part II, Lemma~4.4).
 This representation 
 is supposed to be attached to the unique
 minimal nilpotent coadjoint orbit,
 in the sense that its annihilator in the enveloping algebra $U(\frak g)$
 is  the Joseph ideal if $p+q \ge 8$,
 which is the unique completely prime primitive ideal of
 minimum nonzero Gelfand-Kirillov dimension. 

Our approach based on conformal geometry gives
 a geometric realization of the minimal representation $\varpi^{p,q}$
 for $O(p,q)$.
One of the advantages using  conformal geometry
 is the naturality of the construction
 (see Proposition~2.6),
 which allows us naturally different realizations of $\varpi^{p,q}$,
 not only on the $K$-picture (a compact picture in \S 3),
 but also 
 on the $N$-picture (a flat picture) (see Part III),
 and on the hyperboloid picture (see Part II, \S 7, Corollary~7.2.1),
 together with the Yamabe operator in each realization.
In later sections,
 we shall reduce the branching problems of $\varpi^{p,q}$
 to the analysis on different models
 on which the minimal representation $\varpi^{p,q}$ is realized.

The case of $SO(3,4)$ was treated by \xsab; his method was generalized in
\xtoramin\ to cover all simple groups with admissible minimal orbit, as well
as the case of a local field of characteristic zero.

\def \sc{2}
\sec{}
We write a standard coordinate of $\Bbb R^{p+q}$ as
 $(x, y) = (x_1, \dots, x_p, y_1, \dots, y_q)$.
Let $\Bbb R^{p,q}$ be the pseudo-Riemannian manifold $\Bbb R^{p+q}$
 equipped with the pseudo-Riemannian metric:
$$
     d s^2 = {d x_1}^2 + \cdots + d {x_p}^2 - {d y_{1}}^2 
         - \cdots - d {y_{q}}^2.  
\tag \num.1
$$
We assume $p, q \ge 1$ and define submanifolds of $\Bbb R^{p,q}$ by
$$
\alignat2
     &\Xi  && := \set{(x, y) \in \Bbb R^{p,q}}{|x|=|y|}\setminus \{0\},
\tag \num.2
\\  
      &M  &&:= \set{(x, y) \in \Bbb R^{p,q}}{|x|=|y|=1}
               \simeq S^{p-1} \times S^{q-1}.  
\tag \num.3
\endalignat
$$

We define a diagonal matrix by
 $I_{p,q} := \diag (1, \dots, 1, -1, \dots, -1)$.
The indefinite orthogonal group 
$$
  G = O(p,q) :=\set{g \in GL(p+q, \Bbb R)}{\trans{g} I_{p,q} g = I_{p,q}},
$$ 
 acts isometrically on $\Bbb R^{p,q}$
 by the natural representation,
 denoted by $z \mapsto g\cdot z$ ($g \in G, z \in \Bbb R^{p,q}$).
This action stabilizes the light cone $\Xi$.
The multiplicative group 
 $\Bbb R_+^{\times} := \set{r \in \Bbb R}{r >0}$ acts on $\Xi$
 as a dilation
 and the quotient space $\Xi / \Bbb R_+^{\times}$ is identified with $M$.  
Because the action of $G$ commutes with that of $\Bbb R_+^\times$,
 we can define the action of $G$ 
 on the quotient space $\Xi / \Bbb R_+^{\times}$,
 and also on $M$ through the diffeomorphism
 $M \simeq \Xi / \Bbb R_+ ^{\times}$.
This action will be denoted by
$$
 L_h \: M\to M, x \mapsto L_h x \quad (x \in M, h \in G).
$$ 
In summary,
 we have a $G$-equivariant principal $\Bbb R^\times_+$-bundle:
$$
   \Phi \: \Xi \to M, \ (x,y) \mapsto (\frac{x}{|x|}, \frac{y}{|y|})
                             = \frac{1}{\nu(x,y)} (x,y),
\tag \num.4
$$
 where $\nu \: \Xi  \to \Bbb R_+$ is defined by
$$
     \nu(x,y) = |x| = |y|.  
\tag \num.5
$$
\def \sc{3}
\sec{}
Suppose $N$ is a $p+q-2$-dimensional submanifold of $\Xi$. 
We say $N$ is {\it transversal to rays}
 if $\Phi|_N \: N \to M$ is locally diffeomorphic.
Then,
 the standard pseudo-Riemannian metric on $\Bbb R^{p,q}$
 induces a pseudo-Riemannian metric on $N$
 which has the codimension 2 in $\Bbb R^{p,q}$.
The resulting pseudo-Riemannian metric is denoted by $g_N$,
 which has the signature $(p-1, q-1)$.
In particular,
 $M \simeq S^{p-1} \times S^{q-1}$ itself is transversal to rays,
 and the induced metric $g_{S^{p-1} \times S^{q-1}}
 = g_{S^{p-1}} \oplus (- g_{S^{q-1}})$,
 where $g_{S^{n-1}}$ denotes the standard Riemannian metric on the unit
 sphere $S^{n-1}$.
\proclaim{Lemma~\num}
Assume that $N$ is transversal to rays.
Then $\Phi|_N \: N \to M$ is a conformal map.
Precisely, we have
$$
      (\Phi^* g_M)_z  = \nu(z)^{-2}   (g_{N})_z, 
 \ \text{ for $z = (x,y) \in N$.}
\tag \num.1
$$
\endproclaim
\demo{Proof}
Write the coordinates as $(u_1, \cdots, u_p, v_1, \cdots, v_q)
 =\Phi(x,y) \in S^{p-1} \times S^{q-1}$.
Then  
$$
   \Phi^*\left(d u_j\right) =
   \frac {d x_j} {|x|} -\frac {x_j}{|x|^3} \sum_{i=1}^p x_i d x_i.
$$
Therefore,
 we have
$$
  \align
  \Phi^*\left( \sum_{j=1}^p (d u_j)^2 \right)
   &= |x|^{-2} {\sum_{j=1}^p (d x_j)^2}
      - 2 |x|^{-4} {(\sum_{j=1}^p x_j d x_j)^2}
     + |x|^{-6} {(\sum_{j=1}^p x_j^2)(\sum_{i=1}^p x_i d x_i)^2}
\\
   &= |x|^{-2} {\sum_{j=1}^p (d x_j)^2}
  - |x|^{-4} {(\sum_{j=1}^p x_j d x_j)^2}.  
\intertext{Similarly, we have}
   \Phi^*\left(  \sum_{j=1}^q (d v_j)^2 \right)
    &= |y|^{-2} {\sum_{j=1}^q (d y_j)^2}
                           - |y|^{-4} {(\sum_{j=1}^q y_j d y_j)^2}.  
  \endalign
$$
Because $|x|^2 =|y|^2$
 and $\sum_{j=1}^p x_j d x_j = \sum_{k=1}^q y_k d y_k$,
 we have
$$
     \Phi^*(\sum_{j=1}^p(d u_j)^2 - \sum_{j=1}^q (d v_j)^2)
    =\frac 1 {|x|^2} (\sum_{j=1}^p (d x_j)^2 - \sum_{k=1}^q (d y_k)^2)
$$
Hence, we have proved (\num.1) from
 our definition of
  $g_M$ and $g_N$.
\qed
\enddemo

\def \sc{4}
\sec{}
If we apply Lemma~\ch.3 to the transformation 
 on the pseudo-Riemannian manifold
 $M = S^{p-1} \times S^{q-1}$,
 we have:
\proclaim{Lemma~\num.1}
$G$ acts conformally on $M$.
That is,  
for $h \in G, z \in M$,
 we have 
$$
     L_h^* g_M = \frac 1 {\nu(h \cdot z)^2} g_M
     \quad
     \text{ at $T_z M$}.
$$
\endproclaim
\demo{Proof}
The transformation $L_h\: M \to M$ is 
 the composition of the isometry $M \to h\cdot M, z \mapsto h \cdot z$,
 and the conformal map $\Phi|_{h \cdot M} \: h \cdot M \to M, \xi \mapsto
 \frac{\xi}{\nu(\xi)}$.
Hence Lemma~\num.1 follows.
\qed
\enddemo

Several works in differential geometry treat the connection between the
geometry of a manifold and the structure of its conformal group. 
For the identity
$$
   \Conf(S^{p-1} \times S^{q-1}) = O(p,q),
   \quad
   (p >2, q>2),
$$
see for example \xkobast, Chapter IV.

As in Example~2.2,
 the Yamabe operator on $M = S^{p-1} \times S^{q-1}$ is
 given by the formula:
$$
  \alignat1
     \tilLap{M}
    &=\Delta_{S^{p-1}} - \Delta_{S^{q-1}}
   - \frac{p+q-4}{4(p+q-3)} \left( {(p-1)(p-2)} - (q-1)(q-2) \right)
\\
    &=\left(\Delta_{S^{p-1}}- \frac14 (p-2)^2\right)
     -\left(\Delta_{S^{q-1}}- \frac14 (q-2)^2\right)
\tag \num.1
\\
    &=\left(\tilLap{S^{p-1}}- \frac14\right)
     -\left(\tilLap{S^{q-1}}- \frac14\right).  
  \endalignat
$$
We define a subspace of $C^\infty(S^{p-1}\times S^{q-1})$ by
$$
   \Vpq := \set {f \in C^{\infty} (S^{p-1}\times S^{q-1})}
                {\tilLap{M} f = 0}.
\tag \num.2
$$
By applying Theorem~2.5,
 we have
\proclaim{Theorem~\num.2}
Let $p, q \ge 1$.
For $h \in  O(p,q), z \in M= S^{p-1} \times S^{q-1}$, and $f \in \Vpq$,  
 we define 
$$
    (\varpi^{p,q}(h^{-1})f) (z) := \nu(h \cdot z)^{-\frac {p+q-4} 2}
                                     f(L_h z).  
\tag \num.3
$$
Then $(\varpi^{p,q},\Vpq)$ is a representation of $O(p,q)$.  
\endproclaim

\def \sc{5}
\sec{}
In order to describe the $K$-type formula of $\varpi^{p,q}$,
 we recall the basic fact of spherical harmonics.
Let $p \ge 2$.
The space of spherical harmonics of degree $k \in \Bbb N$
 is defined to be
$$
\alignat1
     \spr{k}{p}
    &= \set{f \in C^{\infty}(S^{p-1})}
          {\Delta_{S^{p-1}} f = - k (k + p -2)f },
\intertext{which is rewritten in terms of 
 $\tilLap{S^{p-1}} = \Delta_{S^{p-1}}-\frac14(p-1)(p-3)$ 
 (see Example~2.2) as}
    &= \set{f \in C^{\infty}(S^{p-1})}
          {\tilLap{S^{p-1}} f
     = \left(\frac{1}{4}- (k + \frac{p -2}2)^2\right) f }.
\tag \num.1
\endalignat
$$
The orthogonal group $O(p)$ acts on $\spr{k}{p}$ irreducibly
 and we have the dimension formula:
$$
   \dim_{\Bbb C} \spr{k}{p} = \pmatrix p +k-2 \\ k \endpmatrix
                    + \pmatrix p + k -3 \\ k-1 \endpmatrix
\tag \num.2
$$
For $p =1$, 
 it is convenient to define representations of $O(1)$ by
$$
    \spr{k}{1} := 
 \cases \Bbb C \ \ \text{(trivial representation)} & \ (k =0)\\
        \Bbb C \ \ \text{(signature representation)} & \ (k =1)\\
         0  & \ (k \ge 2).
 \endcases
$$
Then we have irreducible decompositions as $O(p)$-modules for $p \ge 1$:
$$
     L^2(S^{p-1}) \simeq 
 \overset{\infty\hphantom{M}}\to{\Hsum{k=0\hphantom{M}}}
        \spr{k}{p}
     \quad
     \text{(Hilbert direct sum).}
$$

\def \sc{6}
\sec{}
Here is a basic property of the representation $(\varpi^{p,q},\Vpq)$.  
\proclaim{Theorem~\num.1}
Suppose that $p, q$ are integers with $p \ge 2$ and $q \ge 2$.
\newline
{\rm 1)}\enspace
The underlying \gk-module $(\varpi^{p,q})_K$ of $\varpi^{p,q}$
 has the following $K$-type formula:
$$
    (\varpi^{p,q})_K \simeq \bigoplus 
                     \Sb a, b \in \Bbb N \\ a + \frac{p}2 = b + \frac{q}2
                     \endSb
                      \spr{a}{p} \boxtimes \spr{b}{q}. 
\tag \num.1
$$
\newline
{\rm 2)}\enspace
In the Harish-Chandra parametrization,
 the $\Cal Z(\frak g)$-infinitesimal character of $\varpi^{p,q}$ is given by
$
(1, \frac{p+q}2-2, \frac{p+q}2-3, \dots, 1, 0)
$.
\newline
{\rm 3)}\enspace
$\Vpq$ is non-zero
 if and only if $p + q \in 2 \Bbb N$.
\newline
{\rm 4)}\enspace
If $p + q \in 2 \Bbb N$ and if $(p,q) \neq (2,2)$,
 then $(\varpi^{p,q}, \Vpq)$ is an irreducible representation
 of $G = O(p,q)$
 and the underlying \gk-module
 $(\varpi^{p,q}_K, \Vpq_K)$ is  unitarizable.
\endproclaim
Although Theorem~3.6.1 overlaps with  the results
 of Kostant, Binegar-Zierau, Howe-Tan, Huang-Zhu
 obtained by algebraic methods
 (\xbz, \xhowetan, \xhuzhu, \xkos),
 we shall give a self-contained and new proof from our viewpoint:
 conformal geometry
 and discrete decomposability of the restriction 
 with respect to non-compact subgroups.
The method of finding $K$-types will be generalized
 to the branching law for non-compact subgroups (\S 7, \S 9).
The idea of proving irreducibility seems interesting by its simplicity,
 because we do not need rather complicated computations (cf\. \xbz, \xhowetan)
 but just use the discretely decomposable branching law with 
 respect to $O(p,q') \times O(q'')$.
We shall give two different proofs of the unitarizability
 of $\varpi^{p,q}$ because of the importance of
 \lq\lq small\rq\rq\ representations in the current status of unitary 
 representation theory.
\demo{Proof}
Let $F \in \Vpq \subset C^\infty(M)$.
Then $F$ is developed as
$$
     F = \sum_{a, b \in \Bbb N} F_{a,b}
    \qquad
    (F_{a,b} \in \spr{a}{p} \boxtimes \spr{b}{q}), 
$$
 where the right side converges in the topology of $C^\infty(M)$. 
Applying the Yamabe operator,
 we have
$$
   \tilLap{M}  F 
      = \sum_{a, b \in \Bbb N} 
         \left(-\left(a+\frac{p-2}2\right)^2
        +\left(b+\frac{q-2}{2}\right)^2\right) F_{a,b}.
$$
Since $\tilLap{M}  F = 0$,
 $F_{a,b}$ can be non-zero if and only if
$$
    |a+\frac{p-2}2| = |b+\frac{q-2}{2}|,
\tag \num.2
$$
 whence (1) and (3).
The statement (2) follows from Lemma~\ch.7.2 and (\ch.7.4).
An explicit (unitarizable) inner product for $\varpi^{p,q}$
 will be given in \S \ch.9 (see also Remark~\ch.9, \S 8).

We shall give a simple proof of the irreducibility of $\varpi^{p,q}$
 in Theorem~7.6 by using discretely decomposable
 branching laws to non-compact subgroups (Theorem~4.2 and Theorem~7.1).
\qed
\enddemo

\remark{Remark \num.2}
{\rm 1)}\enspace
$\varpi^{2,2}$ contains the trivial one dimensional representation
 as a subrepresentation.
The quotient $\varpi^{2,2}/\Bbb C$ is irreducible as an $O(2,2)$-module
 and splits into a direct sum of four irreducible $SO_0(2,2)$-modules.
The short exact sequence of $O(2,2)$-modules
 $0 \to \Bbb C \to \varpi^{2,2} \to \varpi^{2,2}/\Bbb C \to 0$
 does not split, 
 and $\varpi^{2,2}$ is not unitarizable as an $O(2,2)$-module.

This case is the only exception that $\varpi^{p,q}$ is not unitarizable
 as a $\Conf(S^{p-1} \times S^{q-1})$-module.
\newline
{\rm 2)}\enspace
The $K$-type formula for the case $p = 1$ or $q=1$
 is obtained by the same method as in Theorem~\num.1.
Then we have 
 that
$$
  \Vpq \simeq \cases
             \Bbb C^4 &\text{ if } (p,q) = (1,1),
\\
   \Bbb C^2 &\text{ if } (p,q) = (1,3), (3,1),
\\
   \{0\}&\text{ if $p =1$ or $q =1$ with $p+q >4$ or if $p+q \notin 2 \Bbb N$}.
              \endcases
$$
 $\Vpq$ consists of locally constant functions
 on $S^{p-1} \times S^{q-1}$ if $(p,q) = (1,1), (1,3)$ and $(3,1)$.
\newline
{\rm 3)}\enspace
In the case of the Kepler problem, i.e. the case of $G = O(4,2)$, the above
$K$-type formula has a nice physical interpretation, namely: the connected
component of $G$ acts irreducibly on the space with positive Fourier
components for the action of the circle $SO(2)$, the so-called {\it positive
energy subspace}; the Fourier parameter $n = 1,2,3,...$ corresponds to the
energy level in the usual labeling of the bound states of the Hydrogen atom,
and the dimension (also called the degeneracy of the energy level) for the
spherical harmonics is $n^2$, as it is in the labeling using angular
momentum and its third component of the wave functions $\psi_{nlm}$.
Here $n$ corresponds to our $b$. 
\endremark

\def\sc{7}
\sec{}
Let us understand $\varpi^{p,q}$ as a subrepresentation 
 of a degenerate principal series.

For $\nu \in \Bbb C$,
 we denote by the space
$$
  S^\nu(\Xi) := \set{f \in C^\infty(\Xi)}{
                  f(t \xi) = t^\nu f(\xi), \ \xi \in \Xi, t > 0}
\tag \num.1
$$
 of smooth functions on $\Xi$ of homogeneous degree $\nu$.
Furthermore,
 for $\epsilon = \pm 1$,
 we put
$$
 S^{\nu, \epsilon}(\Xi) := 
 \set{f \in S^\nu(\Xi)}{f(- \xi) = \epsilon f(\xi), \ \xi \in \Xi}.
$$
Then
 we have a direct sum decomposition
$$
  S^\nu(\Xi) = S^{\nu, 1}(\Xi) + S^{\nu,-1}(\Xi),
$$
 on which $G$ acts by left translations, respectively.
\proclaim{Lemma~\num.1}
The restriction $C^\infty(\Xi) \to C^\infty(M), f \mapsto f|_M$
 induces 
 the isomorphism of $G$-modules between $S^{-\lambda}(\Xi)$ 
 and $(\varpi_\lambda, C^\infty(M))$ (given in (2.5.1))
 for any $\lambda \in \Bbb C$.
\endproclaim
\demo{Proof}
If $f \in S^{-\lambda}(\Xi)$, $h \in G$ and $z \in M$,
 then 
$$
  f(h \cdot z) = f\left(\nu(h \cdot z) \frac{h \cdot z}{\nu(h \cdot z)}\right)
               = \nu(h \cdot z)^{-\lambda} f(L_h z)
               = \left(\varpi_\lambda(h^{-1}) f|_M\right) (z),
$$
 where the last formula follows from the definition (2.5.1) and Lemma~\ch.4.1.
\qed
\enddemo

Let us also identify $S^{\nu, \epsilon}(\Xi)$ with degenerate principal series 
 representations in a standard notation.
The indefinite orthogonal group $G = O(p,q)$
 acts on the light cone $\Xi$ transitively.
We put
$$
   \xi^o := \trans(1, 0, \dots, 0, 0, \dots, 0, 1) \in \Xi.
\tag \num.2
$$
Then the isotropy subgroup at $\xi^o$ is of the form $\Mmax_+ \Nmax$,
 where $\Mmax_+ \simeq O(p-1, q-1)$ and $\Nmax \simeq \Bbb R^{p+q-2}$
 (abelian Lie group).
We set 
$$
  E := E_{1, p+q} + E_{p+q, 1} \in \frak g_0,
$$
 where $E_{i j}$ denotes the matrix unit.
We define an abelian Lie group by
 $\Amax := \exp \Bbb R E  \ (\subset G)$,
 and put
$$
   m_0 := - I_{p+q} \in G.
\tag \num.3
$$
We define $\Mmax$ to be the subgroup generated by $\Mmax_+$ and $m_0$,
 then
$$
   \Pmax := \Mmax \Amax \Nmax
$$
 is a Langlands decomposition
 of a maximal parabolic subgroup $\Pmax$ of $G$.
If $a = \exp(t E)$ $(t \in \Bbb R)$,
 we put $a^\lambda := \exp(t \lambda E)$ for $\lambda \in \Bbb C$.
We put 
$$
   \rho := \frac{p+q-2}2.
$$
For $\epsilon = \pm 1$,
 we define a character  $\chi_\epsilon$ of $\Mmax$ by the composition
$$
  \chi_\epsilon\: \Mmax \to \Mmax/\Mmax_+ \simeq 
  \{1, m_0\} \to \Bbb C^\times,
$$
 such that $\chi_\epsilon(m_0) :=  \epsilon$.
We also write $\operatorname{sgn}$ for $\chi_{-1}$
 and $\boldkey{1}$ for $\chi_{1}$.
We define ${\Cal F}$ to be the ${\Cal A}$, ${\Cal B}$, $C^\infty$ or $\Cal D'$
 valued degenerate principal series by
$$
 \prince{\Cal F}{\lambda}{\epsilon}
 := \set{f \in \Cal F(G)}{f(g m a n) 
 = \chi_\epsilon(m^{-1}) a^{-(\lambda+\rho)} f(g)},
$$
 which has $\Cal Z(\frak g)$-infinitesimal character
$$
     (\lambda,\frac{p+q}2-2,\frac{p+q}2-3,\dots,\frac{p+q}2-[\frac{p+q}2])
\tag \num.4
$$
 in the Harish-Chandra parametrization.
The underlying \gk-module will be denoted by
$
 \princeK{\lambda}{\epsilon}.
$
We note that
  $\princeK{\lambda}{\epsilon}$ is unitarizable 
 if $\lambda \in \sqrt{-1}\Bbb R$.

In view of the commutative diagram of $G$-spaces:
$$
\alignat3
   &G/\Mmax_+ \Nmax &\rarrowsim &\Xi, 
    &&\quad g \Mmax_+ \Nmax \mapsto g \cdot \xi^o.
\\
   &\downarrow &&\downarrow \Phi &&
\\
G/\Pmax \overset{\Bbb Z_2}\to\leftarrow \
  &G/\Mmax_+ \Amax \Nmax &\rarrowsim & M \simeq \Xi/\Bbb R_+^\times
\endalignat
$$
 we have an isomorphism of $G$-modules:
$$
 \prince{C^\infty}{\lambda}{\epsilon}
 \simeq S^{-\lambda -\frac{p+q-2}2, \epsilon}(\Xi).
\tag \num.5
$$

It follows from Theorem~2.5 and Lemma~\num.1 that
 $(\varpi^{p,q}, \Vpq)$ is a subrepresentation of
 $S^{-\frac{p+q-4}2}(\Xi)$.
Furthermore,
 $\varpi^{p,q}(m_0)$ acts on each $K$-type
 $\spr{a}{p}\boxtimes \spr{b}{q}$
 ($a+\frac{p}2 = b + \frac{q}2$) 
 as a scalar
$$
  (-1)^{a+b} = (-1)^{2 a + \frac{p-q}2} = (-1)^{\frac{p-q}2}.
$$
Hence, 
we have the following:
\proclaim{Lemma~\num.2}
 $\varpi^{p,q}$ is a subrepresentation of
$
  S^{a,\epsilon}(\Xi)
$
 with $a =-\frac{p+q-4}2$ and $\epsilon = (-1)^{\frac{p-q}2}$,
 or equivalently,
 of 
$
 \prince{C^\infty}{-1}{(-1)^{\frac{p-q}2}}
$.
\endproclaim
The quotient will be described in (5.5.5).

\remark{Remark \num.3}
1)\enspace
$\varpi^{p,q}$ splits into two irreducible components
 as $SO(p,q)$-modules,
 say $\varpi_\pm^{p,q}$,
 if $p=2$ and $q \ge 4$.
Then,
 $\varpi^{p,q}$ (or $\varpi_\pm^{p,q}$ if $p=2$ and $q \ge 4$)
 coincides with the \lq\lq minimal representations\rq\rq\
 constructed in \xbz, \xkos, \xtoramin.
\newline
{2)}\enspace
In \xbz,
 it was claimed that the  minimal representations of $SO(p,q)$
 are realized in the subspace of
$
  \set{\psi \in C^\infty(S^{p-1} \times S^{q-1})}{\psi(-y) = (-1)^d \psi(y)}
$
 for $d = 2 - \frac{p+q}2$.
But this parity is not correct when both $p$ and $q$ are odd.
\newline{3)}\enspace
Our parametrization of $S^{a,\epsilon}(\Xi)$
 is the same with $S^{a,\epsilon}(X^0)$ in the notation of \xhowetan.
\endremark

\def \sc{8}
\sec{}
Let $p \ge 2$.
The differential operator
$- \Delta_{S^{p-1}} + \frac{(p-2)^2}{4}$ acts on
 the space $\spr{a}{p}$ of spherical harmonics
 as a scalar
$
  a(a+p-2) + \frac14(p-2)^2 = (a +\frac{p-2}2)^2.
$
Therefore,
 we can define a non-negative self-adjoint operator
$$
  D_p
 \: L^2(S^{p-1}) \to L^2(S^{p-1})
\tag \num.1
$$
by
$$
  D_p
 := \left(- \Delta_{S^{p-1}} + \frac{(p-2)^2}{4}\right)^{\frac14}
$$
 with the domain of definition  given by
$$
   \Dom(D_p) :=\set{F = \sum_{a=0}^\infty F_a \in L^2(S^{p-1})}{
 \sum_{a=0}^\infty (a + \frac{p-2}2) \|F_{a}\|^2_{L^2(S^{p-1})}< \infty}.
$$

Here is a convenient criterion which assures a given function to be in 
 $\Dom(D_p)$:
\proclaim{Lemma \num.1}
If $F \in L^2(S^{p-1})$ satisfies $Y F \in L^{2-\frac{2}p}(S^{p-1})$
 for any smooth vector field $Y$ on $S^{p-1}$
  then $F \in \Dom(D_p)$.
Namely, $D_p F$ is well-defined and $D_p F \in L^2(S^{p-1})$.
\endproclaim
In order to prove Lemma~\num.1,
 we recall an inequality due to W\. Beckner:
\proclaim{Fact \num.2 {\rm (\xbec, Theorem~2)}}
Let $1 \le \delta \le 2$ and $F \in L^\delta(S^n)$. 
Let $F = \sum_{k=0}^\infty F_k$
 be the expansion in terms  of spherical harmonics
 $F_k \in \spr{k}{n+1}$,
 which converges in the distribution sense.
Then
$$
   \sum_{k=0}^\infty \gamma_k \| F_k \|^2_{L^2(S^n)}
            \le \| F\|_{L^\delta(S^n)}^2,
\quad
    \gamma_k := \frac{\Gamma(\frac{n}{\delta})\Gamma(k+n-\frac{n}{\delta})}
                 {\Gamma(n-\frac{n}{\delta})\Gamma(k+\frac{n}{\delta})}.
\tag \num.2
$$
\endproclaim
For our purpose,
 we need to give a lower estimate of $\gamma_k$ in Fact~\num.2.
By Stirling's formula for the Gamma function, we have %
$$
  k^{b-a} \frac{\Gamma(k+a)}{\Gamma(k+b)}
  \sim 1 + \frac{(a-b)(a+b-1)}{2k} + \dots
$$
 as $k \to \infty$.
Hence, 
 there exists a positive constant
 $C$ depending only on $n$ and $\delta$ so that
$$
    C k^{n(1-\frac{2}\delta)} \le \gamma_k
\tag \num.3
$$
 for any $k \ge 1$.
Combining (\num.2) and (\num.3), we have:
$$
   C \sum_{k=1}^\infty  k^{n(1-\frac{2}\delta)} \| F_k \|^2_{L^2(S^n)}
            \le \| F\|_{L^\delta(S^n)}^2.
\tag \num.4
$$
Now we are ready to prove Lemma~\num.1.
\demo{Proof of Lemma~\num.1}
Let $\{ X_i\}$ be an orthonormal basis of $\frak o(p)$
 with respect to $(-1) \times$ the Killing form.
The action of $O(p)$ on $S^{p-1}$ induces
 a Lie algebra homomorphism $L\:\frak o(p) \to \frak X(S^{p-1})$.
Then we have $\Delta_{S^{p-1}} = \sum_i L(X_i)^2$.
We write $F = \sum_{k=0}^\infty F_k$
   where $F_k \in \spr{k}{p}$.
We note that $L(X) F_k \in \spr{k}{p}$ for any $k$
 and for any $X \in \frak o(p)$,
 because $\Delta_{S^{p-1}}$ commutes with $L(X)$.
If we apply (\ch.8.4) with $\delta = 2 -\frac{2}{p}$ and $n = p-1$,
 then we have
$$
   C \sum_{k=1}^\infty k^{-1} \| L(X_i) F_k \|^2_{L^2(S^{p-1})}
            \le \| L(X_i) F\|_{L^{2-\frac{2}{p}}(S^{p-1})}^2.
$$
Because $L(X_i)$ is a skew-symmetric operator, we have
$$
   \sum_{i}  \| L(X_i) F_k \|^2_{L^2(S^{p-1})}
  =
   - \sum_{i} (\Delta_{S^{p-1}} F_k, F_k)_{L^2(S^{p-1})}
  = k(k+p-2) \|F_k\|^2_{L^2(S^{p-1})},
$$
and therefore
$$
   C \sum_{k=1}^\infty (k+p-2) \| F_k \|^2_{L^2(S^{p-1})}
            \le \sum_{i} \| L(X_i) F\|_{L^{2-\frac{2}{p}}(S^{p-1})}^2 < \infty.
$$
Hence we have proved that $D_p F$ is well-defined and
$$
 \|D_p F\|^2_{L^2(S^{p-1})} =
 \sum_{k=0}^\infty (k+\frac{p-2}2) \|F_k\|^2_{L^2(S^{p-1})} < \infty.
$$
This completes the proof of Lemma~\num.1
\qed
\enddemo

\def \sc{9}
\sec{}
Let $p \ge 2$ and $q \ge 2$.
We extend $D_p$ to a self-adjoint operator (with the same notation)
 on $L^2(M)$.
Then $D_p$ is a pseudo-differential operator acting on
  $\spr{a}{p} \boxtimes L^2(S^{q-1})$
 as a scalar 
$
 \sqrt{a +\frac{p-2}2}.
$
Likewise, we define $D_q$ as a self-adjoint operator on
 $L^2(S^{q-1})$ and also on $L^2(M)$.
It follows from (\ch.6.2) that 
$$
    D_p  = D_q \quad \text{ on $\Vpq_K$}.
\tag \num.1
$$
Let us consider an explicit unitary inner product for $\varpi^{p,q}$.
We define the Knapp-Stein intertwining operator
$$
A_{\lambda, \epsilon} \:
 \princeK{\lambda}{\epsilon} \to \princeK{-\lambda}{\epsilon}
$$
by
$$
 (A_{\lambda, \epsilon} f)(x)
 := \int_M
 \psi_{\lambda-\rho, \epsilon}([x, b]) \ f(b) \ d b
  \quad
  (x \in M).
$$
Here, $d b$ is the Riemannian measure on $M \simeq S^{p-1} \times S^{q-1}$, $\rho = \frac{
p+q-2}{2}$, and
$$
  \psi_{\nu, \epsilon}(t) := \frac{1}{\Gamma(\frac{2 \nu + 3 - \epsilon}{4})}
                              |t|^\nu \chi_\epsilon(\sgn t).  
$$
Then, 
 $A_{\lambda, \epsilon}$ has an eigenvalue on $\spr{a}{p} \boxtimes \spr{b}{q}$ 
 ($\epsilon = (-1)^{a-b}$)
 as follows (\xkhcrrest):
$$
   \frac{2^{\frac{p+q+2}2-\lambda} \ \pi^{\frac{p+q-1}{2}} 
   (-1)^{[\frac{a-b}2]}
   \ \Gamma(\lambda) \ \Gamma(-B_{\lambda}^{++})}{\Gamma(\frac 1 2) \
       \Gamma(\frac{-2\lambda + p + q -1 -\epsilon}{4}) \
      \Gamma(1+B_\lambda^{--})
     \ \Gamma(1+B_\lambda^{+-})
    \ \Gamma(1+ B_\lambda^{-+})}
\tag \num.2
$$
where for $\epsilon_1, \epsilon_2 = \pm$, we have put
$$
   B_\lambda^{\epsilon_1, \epsilon_2} \equiv 
   B_\lambda^{\epsilon_1, \epsilon_2}(a,b)
 := \lambda - 1 - \epsilon_1(a+\frac{p}2 -1) - \epsilon_2(b + \frac{q}2-1).
\tag \num.3
$$

In particular,
 if $p + q \in 2 \Bbb N$, $p \ge 2$, $q \ge 2$, $(p,q) \neq (2,2)$
 and $\epsilon = (-1)^{\frac{p-q}2}$ ,
 then
$A_{1,\epsilon} \: \princeK{1}{\epsilon} \to \princeK{-1}{\epsilon}$
is non-zero only on
 the submodule $\varpi^{p,q}_K$,
 and has an eigenvalue on $\spr{a}{p} \boxtimes \spr{b}{q}$
 ($a + \frac{p}2 = b + \frac{q}2$) given by
$$
\frac{2^{\frac{p+q}2} \ \pi^{\frac{p+q-1}2} \  
 (-1)^{[\frac{q-p}4]}\ \Gamma(1) \ \Gamma(2a+p-2)}{
\Gamma(\frac{p+q-3-(-1)^{\frac{p-q}2}}4) \ \Gamma(2a+p-1)
 \ \Gamma(1) \ \Gamma(1)}
=
\frac{2^{\frac{p+q-2}2} \ \pi^{\frac{p+q-1}2} 
 (-1)^{[\frac{q-p}4]}}{
\Gamma(\frac{p+q-3-(-1)^{\frac{p-q}2}}4) \ (a+\frac{p}2-1)}
$$
namely,
$$
   \frac{1}{a+\frac{p}2-1} = \frac{1}{b+\frac{q}2-1}
$$
 up to a non-zero constant scalar $(-1)^{[\frac{q-p}4]}c_1$,
 where 
 $$
 c_1 = \frac{\Gamma(\frac{p+q-3-(-1)^{\frac{p-q}2}}4)}{2^{\frac{p+q-2}2} \ \pi^{\frac{p+q-1}2} }.
 \tag \num.4
$$
Because the integration on $G/\Pmax \simeq M/\sim {\Bbb Z}_2$
 gives a $G$-invariant non-degenerate sesquilinear form 
$$
\princeK{1}{\epsilon} \times \princeK{-1}{\epsilon} \to \Bbb C,
$$
 the underlying $(\frak g,K)$-module $\Vpq_K$
 ($p, q \ge 2$, $(p,q) \neq (2,2)$, $p+q \in 2 \Bbb N$)
 is unitarizable
 with the inner product
$$
  (f_1, f_2) := \int_M (a+\frac{p}2-1) f_1 \overline{f_2} \ d \omega,
   \quad
   f_1, f_2 \in \spr{a}{p} \boxtimes \spr{b}{q}, %
$$
namely,
$$
  (f_1, f_2) = \int_M (D_p f_1) \overline{D_p f_2} \ d \omega
              =  \int_M (D_q f_1) \overline{D_q f_2} \ d \omega,
   \quad
   f_1, f_2 \in \Vpq_K
\tag \num.5
$$
 where $d \omega$ is the standard measure on $M$.

We denote by $\overline\Vpq$ the Hilbert completion of $\Vpq$,
 on which $G$ acts as an (irreducible) unitary representation of $G$.
We shall use the same notation $\varpi^{p,q}$ to denote this unitary
 representation.

In view of \S \ch.8,
 we can describe $\overline\Vpq$ as follows:
Let $\Cal V$ be the Hilbert space of the completion of
 $C^\infty(M)$ by the norm defined by
$$
  \| F\|^2_{L^2(M)} + \|(D_p+D_q) F\|^2_{L^2(M)} 
   \quad
   \text{ for }  
   F \in C^\infty(M).
$$
Then,
 $\Cal V$ is a dense subspace of $L^2(M)$ and
 $\Cal V = \Dom(D_p) \cap \Dom(D_q)$.
With this notation,
 the closure $\overline{\Vpq}$ is characterized directly by the following: 
\proclaim{Theorem ~\num} 
$$
  \overline\Vpq 
 = \set{f \in \Cal V}{D_p f = D_q f}
 = \set{f \in \Cal V}{\tilLap{M} f =0},
$$
 where $\tilLap{M} f=0$ is in the distribution sense.
The norm for the unitary representation $\overline\Vpq$
 equals
$
\frac14 \|(D_p+D_q) F\|^2_{L^2(M)}. 
$
Namely,
 if
 $F = \sum_{a} F_{a,b} \in \overline\Vpq$
 with $F_{a,b} \in \spr{a}{p}
 \boxtimes \spr{b}{q}$ and $b = a + \frac{p-q}2$,
 then $D_p F = \sum_{a} \sqrt{a + \frac{p-2}2} F_{a,b}$ and
$$
   \| F \|^2 = \sum_{a} (a + \frac{p-2}2) \|F_{a,b}\|^2_{L^2(M)}.
$$
Here the sum is taken over integers $a$ with $a \ge \max(0,\frac{p-q}2)$.
\endproclaim
 
\remark{Remark}
By comparing the construction of \xbz, 
 (\num.4) coincides with the formula obtained by Binegar-Zierau 
 by a different method (see Remark~\ch.7.3~(1)).
These authors defined a similar operator $\Cal D_n$ of \xbz, page 249;
we remark that
 there is a typographical error
 in the definition of $\Cal D_n$; 
 $n-2$ should read $(n-2)^2$.
Then, the square $D_p^2$ of our operator corresponds to
 $\Cal D_p$ if $p \neq 2$; $|\Cal D_p|$ if $p=2$, 
 with the notation in \xbz.
\endremark

\def \sc{10}
\sec{}
The following lemma is rather weak, 
 but is clear from the Sobolev estimate.
\proclaim{Lemma~\num}
Suppose $W$ is an open set of $M$ such that the measure
 of $M \setminus W$ is zero.
Suppose $F$ is a $C^\infty$ function on $W$
 satisfying $\tilLap{M} F = 0$ on $W$.
If $F \in L^2(M)$ and if $Y Y' F \in L^{1}(M)$ for any $Y, Y' \in \frak X(M)$
 (differentiation in the sense of the Schwartz distributions),
 then $F \in \overline{\Vpq}$.
\endproclaim

\redefine\cite{{}}

\enddocument